\documentclass[
    ,final            
  ]
{aipproc}
\usepackage{amsmath, amssymb, bm}
\usepackage{tikz}
\layoutstyle{8x11single}

\def\pmatrix{\left(\begin{array}}
\def\endpmatrix{\end{array}\right)}

\def\wmQ{\widetilde{\mathcal{Q}}}
\def\wmA{\widetilde{\mathcal{A}}}
\def\sech{\mathrm{sech}}
\DeclareMathAlphabet{\mathcal}{OMS}{cmsy}{m}{n}
\SetMathAlphabet{\mathcal}{bold}{OMS}{cmsy}{b}{n}

\begin{document}

\title{Bespoke finite difference methods that preserve two local conservation laws of the modified KdV equation.}


\author{Gianluca Frasca-Caccia}{
  address={School of Mathematics Statistics and Actuarial Science (SMSAS), University of Kent, UK}
}

\begin{abstract} 
By exploiting the fact that conservation laws form the kernel of a discrete Euler operator, we use a recently introduced symbolic-numeric approach to construct a new class of finite difference methods for the modified Korteweg-de Vries (mKdV) equation, that preserve the local conservation laws of mass and energy.
\end{abstract}

\maketitle


\section{INTRODUCTION} In the last decades there has been a growing interest in the development of structure-preserving algorithms for the numerical solution of partial differential equations (PDEs). Particularly, for the solution of Hamiltonian PDEs, multisymplectic integrators, able to preserve a discrete conservation law of symplecticity, have been widely considered \cite{AschMcLac2,AschMcLac,AK10,Bridges,BridgesReich,LeimkuhlerReich}. Another popular strategy is to use a method of lines approach to obtain a semidiscrete Hamiltonian systems of ODEs which can be integrated (in time) by a symplectic \cite{AschMcLac,BridgesReich,LeimkuhlerReich} or an energy-conserving method \cite{BBFI18,BFI15,BI,BI18}. 

In this paper, we consider a different procedure to construct bespoke finite difference methods that preserve multiple local conservation laws of a given PDE. This approach, first introduced in \cite{GrantHydon} and further developed in \cite{FH,Grant}, has three main advantages. First, it does not require the PDE to have any special structure, so, for instance, it can be applied also to non-Hamiltonian PDEs (see \cite{FH}). Second, preserving {\em local} conservation laws implies the preservation of the corresponding global invariants (given suitable boundary conditions) and gives a stricter constraint in general, as the converse is not true. Finally, it can be used to seek methods that preserve any number of conservation laws. Nevertheless, this in general increases the complexity of the method.

We present this approach by considering a PDE with two independent variables (but it can be generalized),
\begin{equation}\label{PDE}
\mathcal{A}(x,t,[u])=0,
\end{equation}
where $[u]$ denotes the dependent variable $u(x,t)$ and finitely many of its derivatives. 

A conservation law of (\ref{PDE}) is a divergence expression
\begin{equation}\label{CLAW}
\mbox{Div\,} \mathbf{F}\equiv D_x\{F(x,t,[u])\}+D_t\{G(x,t,[u])\},
\end{equation}
which is zero on all solutions of (\ref{PDE}). Here $D_x$ and $D_t$ denote the total derivatives with respect to $x$ and $t$ respectively. The components $F$ and $G$ are referred to as the {\em flux} and {\em density} respectively.

If the conservation law (\ref{CLAW}) amounts to
$$\mbox{Div\,} \mathbf{F}=\mathcal{QA},$$
it is said to be in {\em characteristic form} and the multiplier function $\mathcal{Q}(x,t,[u])$ is called a {\em characteristic} of the conservation law.

In order to discretize the PDE (\ref{PDE}) we introduce a uniform mesh. With respect to a generic lattice point $\mathbf{n}=(m,n),$ the grid points are 
$x_i=x(m)+i\Delta x,$ and $t_j=t(n)+j\Delta t$. We denote with $u_{i,j}$ the approximation of $u(x_i,t_j)$, and we introduce the forward shift operators $S_m$ and $S_n$ defined as:
$$S_m:(m,n,x_i,t_j,u_{i,j})\mapsto (m+1,n,x_{i+1},t_j,u_{i+1,j}),\quad S_n:(m,n,x_i,t_j,u_{i,j})\mapsto (m,n+1,x_{i},t_{j+1},u_{i,j+1}).$$
Denoting with $I$ the identity operator, the forward difference, $D_m$, $D_n$, and the forward average, $\mu_m$, $\mu_n$, are defined by
\[ D_m=\tfrac{1}{\Delta x}(S_m-I),\quad D_n=\tfrac{1}{\Delta t}(S_n-I),\quad \mu_m=\tfrac{1}{2}(S_m+I),\quad \mu_n=\tfrac{1}{2}(S_n+I).
\]
Approximating the derivatives in (\ref{PDE}) by means of suitable finite differences, yields a partial difference equation (P$\Delta$E),
\begin{equation}\label{PDeltaE}
\widetilde{\mathcal{A}}(x_i,t_j,u_{i,j})=0.
\end{equation}
Here and henceforth, tildes represent discretizations of the corresponding continuous terms. We seek schemes having the following discrete analogue of each conservation law:
\begin{equation}\label{DCLAW}
\mbox{Div\,} \widetilde{\mathbf{F}}\equiv D_m\{\widetilde F(x_i,t_j,u_{i,j})\}+D_n\{\widetilde G(x_i,t_j,u_{i,j})\},
\end{equation}
which is zero on all solutions of (\ref{PDeltaE}). The functions $\widetilde{F}$ and $\widetilde{G}$ are respectively the flux and the density of the conservation law (\ref{DCLAW}). 
A conservation law (\ref{DCLAW}) is in characteristic form if
\[\mbox{Div\,} \widetilde{\mathbf{F}}=\widetilde{\mathcal{Q}}\widetilde{\mathcal{A}},\]
where the multiplier function $\wmQ(x_i,t_j,u_{i,j})$ is called the {\em characteristic} \cite{Hydonbook}.

A crucial result for our purposes is due to Kuperschmidt \cite{Kuperschmidt} and states that the set of all discrete divergence expressions (\ref{DCLAW}) is the kernel of the difference Euler operator,
\[\mathbf{E}=\sum_{i,j}S_m^{-i}S_n^{-j}\frac{\partial}{\partial u_{i,j}},\]
(see also \cite{HydonMans} for the generalisation of this result).
Consequently, if $\widetilde{\mathcal{Q}}\approx\mathcal{Q}$ is such that $\mathbf{E}(\wmQ\wmA)=0$, there exists $\widetilde{\mathbf{F}}\approx\mathbf{F}$ such that $\wmQ\wmA=\mbox{Div\,} \widetilde{\mathbf{F}}$, corresponding to a conservation law of which $\wmQ$ is the characteristic.

The approach in \cite{FH,Grant,GrantHydon} to find schemes that preserve conservation laws is straightforward. Choose a stencil of points and consider the most general discretizations  $\wmA$ of the PDE and $\wmQ$ of the characteristic of the desired conservation law. These discretizations depend on a number of free parameters. In order to preserve the conservation law, set the parameters in such a way to satisfy the constraints given by the algebraic condition $\mathbf{E}(\wmQ\wmA)=0$. If not all the free parameters have been set, this procedure can be iterated for multiple characteristic $\wmQ_\ell$, to obtain the preservation of all the corresponding conservation laws, provided that the system of algebraic equations obtained can be solved. Finally, consistency conditions are applied to ensure that $\wmA$ converges to $\mathcal{A}$ and each $\wmQ_\ell$ converges to $\mathcal{Q}_\ell$ as the stepsizes tend to zero, giving further constraints on the free parameters.

A crucial step of this procedure is to solve, for each characteristic $\wmQ_\ell$, the condition
\begin{equation}\label{E(CA)}
\mathbf{E}(\wmQ_\ell\wmA)=0.
\end{equation}
This has to be done symbolically and it is not easy in general. In particular, if the non-linearity in $\mathcal{A}$ and $\mathcal{Q}_\ell$ is of polynomial type, (\ref{E(CA)}) amounts to a large system of polynomial equations. This, in principle, can be solved by finding a Groebner basis, but the calculation may take a huge amount of memory and a very long computation time even when the nonlinearity is only quadratic and using a compact stencil. 

Nevertheless, the complexity of the symbolic calculations may be reduced by restricting the dependence of few key quantities in $\wmA$ and $\wmQ_\ell$, on only points of the most compact sub-stencil that allows the desired order of accuracy. In particular, approximating nonlinear terms using as few points as possible, may be effective to the point of being able to solve (\ref{E(CA)}) with a fast symbolic computation that does not need a Groebner basis (see \cite{FH}).

In the next section we use the procedure described above to find bespoke finite difference schemes preserving two conservation laws of a Hamiltonian PDE presenting a cubic nonlinearity.

\section{MODIFIED KDV EQUATION}
In this section we use the procedure described above to develop conservative finite difference schemes for the modified Korteweg-de Vries (mKdV) equation
\begin{equation}\label{mKdV}
\mathcal{A}\equiv u_t+u^2u_x+u_{xxx}=0, \quad (x,t)\in\Omega\equiv[a,b]\times[0,\infty).
\end{equation}
Equation (\ref{mKdV}) has infinitely many conservation laws. The first three are
\begin{align}\label{CL1}
&D_t(G_1)+D_x(F_1)\equiv D_t(u)+D_x\left(\tfrac{1}3u^3+u_{xx}\right)=0,\\\label{CL2}
&D_t(G_2)+D_x(F_2)\equiv D_t(\tfrac{1}2 u^2)+D_x\left(\tfrac{1}4u^4+uu_{xx}-\tfrac{1}2u_x^2\right)=0,\\\label{CL3}
&D_t(G_3)+D_x(F_3)\equiv D_t(\tfrac{1}{12} u^4+\tfrac{1}2uu_{xx})+D_x\left(\tfrac{1}2\left(\tfrac{1}3u^3+u_{xx}\right)^2+u_xu_t-uu_{xt}\right)=0.
\end{align}
These are known as the local conservation laws of mass, momentum and energy, respectively and, assuming conservative boundary conditions, integration in space yields the conservation of three global invariants. The conservation laws (\ref{CL1})--(\ref{CL3}) can be written in characteristic form with characteristics, respectively, 
\[\mathcal{Q}_1=1,\quad \mathcal{Q}_2=u, \quad \mathcal{Q}_3=\tfrac{1}3u^3+u_{xx}.\]

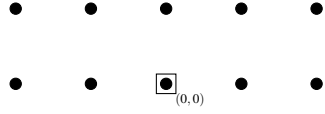
\begin{figure}
\begin{tikzpicture}
\draw[fill] (0,0) circle (0.075);
\draw[fill] (1,0) circle (0.075);
\draw[fill] (2,0) circle (0.075);
\draw[fill] (3,0) circle (0.075);
\draw[fill] (4,0) circle (0.075);
\draw[fill] (0,1) circle (0.075);
\draw[fill] (1,1) circle (0.075);
\draw[fill] (2,1) circle (0.075);
\draw[fill] (3,1) circle (0.075);
\draw[fill] (4,1) circle (0.075);
\draw (1.87,-0.13) rectangle (2.13,0.13);
\node [below right] at (2,0) {\tiny{$(0,0)$}};

\end{tikzpicture}
\caption{10-point rectangular stencil.}
\label{stencil10}
\end{figure}
Our purpose here is to find second-order accurate finite difference schemes defined on the 10-point stencil in Figure~\ref{stencil10} that preserve both the conservation laws (\ref{CL1}) and (\ref{CL3}). Henceforth grid points are labelled with respect to the lattice point denoted with a square in Figure~\ref{stencil10}. Setting $\wmQ_1=1$, the second-order approximations of the conservation laws (\ref{CL1}) and (\ref{CL3}) are taken to be of the form
\begin{align*}
\wmA&=D_m(\widetilde{F_1})+D_n(\widetilde{G_1})=0,\qquad \wmQ_3\wmA=0.
\end{align*}
Since $\wmA$ is defined to be a discrete conservation law, $\mathbf{E}(\wmQ_1\wmA)=\mathbf{E}(\wmA)=0,$ for any discretizations $\widetilde{F_1}$ and $\widetilde{G_1}$. Nevertheless, as both $\wmA$ and $\wmQ_3$ include a cubic nonlinearity, the solution of $\mathbf{E}(\wmQ_3\wmA)=0$ cannot be easily tackled without introducing some simplifying assumption on the approximations (see \cite{FH}). Approximating both $G_1$ and the cubic term in $\mathcal{Q}_3$ on compact sub-stencils makes the symbolic computations fast. Taking into account that the approximations have to be second-order, this implies that $\widetilde{G_1}$ depends only on the grid point $(0,0)$ and the approximation of the cubic term in $\mathcal{Q}_3$ depends only on the two lattice points $(0,0)$ and $(0,1)$. Hence we consider here discretizations of $\mathcal{Q}_3, {G_1}$ and ${F_1}$ of the form
\begin{align*}\nonumber 
\wmQ_3&=\alpha_1u_{0,0}^3+\alpha_2u_{0,1}u_{0,0}^2+\alpha_3u_{0,0}u_{0,1}^2+\alpha_4u_{0,1}^3+\frac{1}{\Delta x^2}\sum_{i=-2}^2\sum_{j=0}^1\beta_{i,j}u_{i,j},\qquad \widetilde{G_1}=u_{0,0},\\ 
\widetilde{F_1}&=\sum_{i=-2}^1\left(\sum_{k=i}^1\sum_{r=k}^1\sum_{j=0}^1\gamma_{i,j,k,j,r,j}u_{i,j}u_{k,j}u_{r,j}+\sum_{k=-2}^1\sum_{r=k}^1\left(\gamma_{i,1,k,0,r,0}u_{i,1}u_{k,0}u_{r,0}+\gamma_{i,0,k,1,r,1}u_{i,0}u_{k,1}u_{r,1}\right)+\frac{1}{\Delta x^2}\sum_{j=0}^1\xi_{i,j}u_{i,j}\right)
.\end{align*}

We find the undetermined coefficients $\alpha,\,\beta,\,\gamma,$ and $\xi$ by symbolically solving $\mathbf{E}(\wmQ_3\wmA)=0$ and imposing the consistency conditions giving second-order accuracy. This yields a one-parameter family of second-order schemes,
\begin{align*}
&\mbox{EC}(\lambda)\equiv D_m(\widetilde{F_1})+D_n(\widetilde{G_1})=0,\qquad \widetilde{F_1}=\mu_m\varphi_{-1,0},\qquad \widetilde{G_1}=u_{0,0},\\
& \varphi_{-1,0}=\tfrac{1}3(\mu_n u_{-1,0}^2)(\mu_n u_{-1,0})+\mu_nD_m^2 u_{-2,0}+\lambda\Delta x^2D_mD_n\mu_mu_{-2,0}.
\end{align*}
Clearly, each of these schemes preserves the local conservation law of the mass. The numerical solution also satisfies the following discrete energy conservation law
\begin{align*}
&\wmQ_3\wmA=D_m(\widetilde{F_3})+D_n(\widetilde{G_3})=0,\qquad \wmQ_3=\varphi_{0,0},\qquad \widetilde{G_3}=\tfrac{1}{12}u_{0,0}^4+\tfrac{1}2 u_{0,0}(D_m^2u_{-1,0})\\
&\widetilde{F_3}=\tfrac{1}2\left(\varphi_{-1,0}\varphi_{0,0}+(D_m\mu_nu_{-1,0})(D_n\mu_mu_{-1,0})-(\mu_m\mu_nu_{-1,0})(D_mD_nu_{-1,0})+\lambda\Delta x^2(D_nu_{0,0})(D_nu_{-1,0})\right),
\end{align*}
where $\widetilde{F_3}$ and $\widetilde{G_3}$ have been reconstructed from the characteristic (see \cite{Hydonbook}).
As the free parameter $\lambda$ appears as a factor of a second order perturbation in each scheme, one may be able to find an optimal value that reduces the local truncation error. However, this depends on the particular problem and no choice of $\lambda$ gives a higher order method.

\section{NUMERICAL TEST}
In this section we present a numerical test that shows the effectiveness of the schemes developed in the previous section compared with two well-known methods both satisfying only a discrete local conservation law of the mass. These are the multisymplectic and the narrow box schemes proposed in \cite{AK10} and inspired to the analogous methods developed by Ascher and McLachlan in \cite{AschMcLac2,AschMcLac} for the KdV equation. These schemes are both implicit and defined on an 8-point stencil, which is the most compact for discretizing the mKdV equation. We consider equation (\ref{mKdV}) setting periodic boundary conditions over $[-20,20]$ for $t\in [0,2]$ and initial condition given by the single-soliton solution on~$\mathbb{R}$,
\[u(x,t)=\sqrt{30}\,\sech\left(\sqrt{5}x-5\sqrt{5}(t-1)\right).\]
Each scheme is solved with stepsizes $\Delta x=0.1$ and $\Delta t=0.01$. In Table~\ref{results}, we compare the relative error in the solution at the final time, and the error in the conservation laws (\ref{CL1})-(\ref{CL3}) measured by the maximum value of the absolute error in the global invariants at every iteration step. These are denoted by $\mbox{Err}_1$, $\mbox{Err}_2$ and $\mbox{Err}_3$, respectively. 
We show the results for two particular EC$(\lambda)$ schemes corresponding to the values of the parameter $\lambda=0.023$ and $\lambda=-0.07$ minimizing the error in the solution and in the global momentum, respectively. EC$(0.023)$ is by far the most accurate scheme. Nevertheless, choosing the free parameter to minimize the non-preserved conservation law, yields a solution error about 2.5 and 5.5 times smaller than the narrow box and the multisymplectic scheme, respectively.

\quad\\
{\em Aknowledgement.} The author would like to express his gratitude to prof. Peter E. Hydon (University of Kent) for valuable discussions and support throughout this project.

\begin{table}[t]
\caption{Maximum value of the absolute error in global invariants at every step and relative error in the solution at the final time.}
\label{results}
\small
\begin{tabular}{ccccc}
\hline
Method &  $\text{Err}_1$ & $\text{Err}_2$ & $\text{Err}_3$  &  Error in solution\\ 
\hline
$\mbox{EC}(0.023)$	& 1.69e-14  & 1.41e-04 & 1.24e-13 &  0.0036\\
$\mbox{EC}(-0.07)$	& 1.78e-14  & 9.50e-05 & 1.88e-13 &  0.0587\\
Multisymplectic & 3.73e-14  & 6.56e-04 & 0.7423 &  0.3225\\
Narrow box  & 2.22e-14  & 0.0013 & 0.7836 &  0.1443\\
\hline
\end{tabular}
\end{table}

\bibliographystyle{aipproc}

\end{document}